\documentclass[12pt,a4paper]{article}
\usepackage{amsmath,amsthm,amssymb, mathrsfs, tocloft}
\usepackage{hyperref}
\input xy
\xyoption{all}
\usepackage{epsfig}

\newtheorem{thm}{Theorem}[section]

 \theoremstyle{definition}

\theoremstyle{remark}
\newtheorem{remark}[thm]{Remark}

\numberwithin{equation}{section}

\newcommand{\ben}{\begin{enumerate}}

\newcommand{\een}{\end{enumerate}}
\newcommand{\bit}{\begin{itemize}}
\newcommand{\eit}{\end{itemize}}

\begin{document}

\title
{Recovering Functions from the Spherical Mean Transform with Data on an Ellipse Using Eigenfunction Expansion in Elliptical Coordinates}
\author{Yehonatan Salman \\ Weizmann Institute of Science \\ Email: salman.yehonatan@gmail.com}
\date{}
\maketitle

\begin{abstract}

The aim of this paper is to introduce a new inversion procedure for recovering functions, defined on $\Bbb R^{2}$, from the spherical mean transform, which integrates functions on a prescribed family $\Lambda$ of circles, where $\Lambda$ consists of circles whose centers belong to a given ellipse $E$ on the plane. The method presented here follows the same procedure which was used by S. J. Norton in \cite{22} for recovering functions in case where $\Lambda$ consists of circles with centers on a circle. However, at some point we will have to modify the method in \cite{22} by using expansion in elliptical coordinates, rather than spherical coordinates, in order to solve the more generalized elliptical case. We will rely on a recent result obtained by H.S. Cohl and H.Volkmer in \cite{8} for the eigenfunction expansion of the Bessel function in elliptical coordinates.

\end{abstract}

\section{Introduction and Mathematical Background}

\hskip0.6cm The main object of study in this paper is the spherical mean transform (SMT for short) which integrates functions on circles in the plane. For a function $f$, defined on $\Bbb R^{2}$, we denote its SMT by $Rf$ which is defined explicitly as follows
\vskip-0.2cm
$$\hskip-8cm Rf:\Bbb R^{2}\times\Bbb R^{+}\rightarrow\Bbb R,$$
$$\hskip-3cm(Rf)(x,r) = \int_{\Bbb S^{1}}f(x + ry)rdS(y), x\in\Bbb R^{2}, r\in\Bbb R^{+}.$$
Here, $\Bbb R^{+}$ is the ray $[0,\infty)$, $\Bbb S^{1}$ is the unit circle in $\Bbb R^{2}$, i.e., $\Bbb S^{1} = \{x\in\Bbb R^{2}:|x| = 1\}$ and $dS(y)$ is the standard infinitesimal length measure on $\Bbb S^{1}$.

The SMT has been found out to be a useful tool in many different disciplines in applied mathematics and science such as thermo
and photoacoustic tomography (\cite{3, 15, 17, 25, 26, 27}), reflectivity tomography (\cite{18}), geophysics (\cite{4,9}), radar and sonar imaging (\cite{6, 7, 19}) and approximation theory (\cite{1}).

One of the most important inverse problems regarding the SMT is the determination of a function $f$ from its SMT $Rf$ in case where $Rf$ is restricted to a cylindrical surface of the from $\Gamma\times\Bbb R^{+}$ where $\Gamma$ is a closed curve in $\Bbb R^{2}$. That is, we assume that our data consists of the integrals of $f$ on any circle with a center on the curve $\Gamma$ while no restriction is imposed on the set of radii and where our goal is to reconstruct the function $f$. The above inverse problem was deeply investigated in the last few decades since it serves as the base model for the wave propagation in thermo and photoacoustic tomography. This wave propagation occurs when a short duration electromagnetic pulse intersects an object $\Omega$ and which is then can be measured by detectors which are scattered around $\Omega$ (see \cite{3, 27} and references therein).

In the last few decades, explicit inversion methods for the case where $\Gamma$ is a circle, or more generally an ellipse (or ever more generally an ellipsoid in $\Bbb R^{n}$), where derived in \cite{2, 10, 11, 12, 13, 14, 16, 20, 21, 22, 23, 24, 28}. The inversion methods obtained in \cite{2, 10, 11, 16, 21, 22, 28} for the case where $\Gamma$ is a circle (or more generally a hypersphere) can be roughly divided into three different types.

\begin{itemize}

\item Inversions formulas of back-projection type where in order to recover a function $f$, at a certain point $x$, one averages the values of $Rf$ which is now restricted to the set of all circles (or more generally hyperspheres) with centers on $\Gamma$ and which pass through $x$ (for example \cite{21, 28}).

\item Inversion formulas that use a change of order of integration in order to obtain from the SMT $Rf$ of $f$ an integral equation of the form $G(x) = \int f(y)L(x - y)dy$ where $L$ is a fundamental solution for the Laplace operator in $\Bbb R^{2}$ (or more generally in $\Bbb R^{n}$). Then, $f$ can be recovered by taking the Laplace operator on $G$ (for example \cite{2, 10, 11}).

\item Inversion formulas obtained by expanding $f$ and $Rf$ into spherical harmonics in order to reduce the inversion problem for each term in the expansions of $f$ and $Rf$ respectively and then using each term of $Rf$ in order to recover the corresponding spherical harmonic term of $f$. This method was used by S. J. Norton in \cite{22}.

\end{itemize}

The two first inversion methods, which use back-projection procedures and fundamental solutions for the Laplace operator, have recently been generalized in \cite{12, 13, 14, 20, 24} for the case where $\Gamma$ is an ellipse, or more generally an ellipsoid in $\Bbb R^{n}$, with only slight modifications. However, the inversion procedure obtained by Norton in \cite{22}, which uses expansion into spherical harmonics, has never been modified to include the case where $\Gamma$ is an ellipse. The main problem in generalizing the method in \cite{22} is that when using spherical harmonics expansion one should assume, in order to extract each term in the expansion, an invariance of our collection of data with respect to the group of rotations corresponding to our space.

In the case where our family $\Lambda$ of circles of integration consists of circles with centers on the curve $\Gamma$, where $\Gamma$ itself is a circle, then obviously $\Lambda$ is invariant with respect to the group of rotations with respect to the center of $\Gamma$. However, when $\Gamma$ is an arbitrary ellipse then $\Lambda$ does not satisfy this invariance property anymore. Hence, for the case where $\Gamma$ is an ellipse we will modify the method used by Norton in \cite{22} by using expansions into eigenfunctions of the Laplace operator in elliptical coordinates. For this, we will use the recent result obtained by H.S. Cohl and H.Volkmer in \cite{8} for the expansion into eigenfunctions in elliptical coordinates of $J_{0}(|x|)$ where $J_{0}$ denotes the Bessel function of order zero.

It should be mentioned that since, in this paper, expansion into eigenfunctions is used then our inversion formula will be given as an infinite sum rather than in a closed form and thus it has more of a theoretical importance rather than practical. However, the main emphasize of this paper is to show how expansion into eigenfunctions in elliptical coordinates is a proper tool for solving problems in Integral Geometry with the vision that similar problems could be solved with the same approach.\\

Before stating the main result of this paper we will introduce some basic notations and definitions.

Let us fix a point $\xi_{0}$ in $\Bbb R^{+}$ and denote by $E$ the following ellipse
\begin{equation}\hskip-4cm E = \left\{(x_{1},x_{2})\in\Bbb R^{2}:\frac{x_{1}^{2}}{\cosh^{2}\xi_{0}} + \frac{x_{2}^{2}}{\sinh^{2}\xi_{0}} = 1\right\}\end{equation}
in $\Bbb R^{2}$. For a point $x\in\Bbb R^{2}$ denote by $x(\xi,\eta)$ its presentation in elliptical coordinates, that is
$$x(\xi, \eta) = (\cosh(\xi)\cos(\eta), \sinh(\xi)\sin(\eta)), \xi\in\Bbb R^{+}, \eta\in[-\pi,\pi),$$
where the quantities $\xi = \xi(x)$ and $\eta = \eta(x)$ are uniquely determined by $x$. In elliptical coordinates the Helmholtz equation
\vskip-0.2cm
$$\hskip-7cm(\triangle_{x} + k^{2})U(x) = 0$$
is given by
\vskip-0.2cm
$$\hskip-2.5cm\left(\frac{\partial^{2}}{\partial\xi^{2}} + \frac{\partial^{2}}{\partial\eta^{2}} + k^{2}(\cosh^{2}\xi - \cos^{2}\eta)\right)u(\xi, \eta) = 0$$
where $u(\xi, \eta) = U(x).$ Using separation of variables $u(\xi, \eta) = u_{1}(\xi)u_{2}(\eta)$ yields the following two equations
\vskip-0.2cm
$$\hskip-4cm - u_{1}^{''}(\xi) + (\sigma - 2q\cosh(2\xi))u_{1}(\xi) = 0,$$
\begin{equation}
\hskip-4cm u_{2}^{''}(\eta) + (\sigma - 2q\cos(2\eta))u_{2}(\eta) = 0,
\end{equation}
where $q = \frac{1}{2}k^{2}$. In \cite{8} it was proved that the Mathieu's equation (1.2) has a solution if and only if $\sigma$ is one of its eigenvalues $a_{n}(q), n = 0,1,2,...$ which correspond to even $2\pi$ periodic eigenfunctions $\textrm{ce}_{n}(\eta, q)$ or if $\sigma$ is one of its eigenvalues $b_{n}(q), n = 1,2,...$ which correspond to odd $2\pi$ periodic eigenfunctions $\textrm{se}_{n}(\eta, q)$. For notation convenience we will also define $\textrm{se}_{0}(\eta, q) = 0$.

The family
\vskip-0.2cm
$$\hskip-4.5cm\Theta = \{\textrm{ce}_{n}(.,q)\}_{n = 0}^{\infty}\cup\{\textrm{se}_{n}(.,q)\}_{n = 0}^{\infty}$$
is an orthogonal set of functions (see \cite[page 654]{5}) on the interval $\eta\in[-\pi,\pi)$ which are uniquely determined by the following equations
\begin{equation}
\hskip-4cm \int_{-\pi}^{\pi}\textrm{ce}_{n}^{2}(\eta, q)d\eta = \int_{-\pi}^{\pi}\textrm{se}_{n}^{2}(\eta, q)d\eta = \pi.
\end{equation}
Denote by $\rho(\xi,\eta,\lambda,\theta) = |x(\xi, \eta) - x(\lambda, \theta)|$ the distance between the points $x(\xi,\eta)$ and $x(\lambda,\theta)$. The following expansion
\vskip-0.2cm
$$\hskip-1.5cm J_{0}(k\cdot\rho(\xi, \eta, \lambda, \theta)) = \frac{1}{\pi}\sum_{n = 0}^{\infty}\left[\mu_{n}(q)\textrm{ce}_{n}(i\xi, q) \textrm{ce}_{n}(\eta,q)\textrm{ce}_{n}(i\lambda, q)\textrm{ce}_{n}(\theta,q)\right.$$
\begin{equation}
\left. \hskip3.3cm + \upsilon_{n}(q)\textrm{se}_{n}(i\xi, q)\textrm{se}_{n}(\eta,q)\textrm{se}_{n}(i\lambda, q)\textrm{se}_{n}(\theta,q)\right]
\end{equation}
for the Bessel function $J_{0}$ was proved in \cite[Theorem 4.2]{8} where $q = \frac{1}{2}k^{2}$ and where the functions $\mu_{n}(q)$ and $\upsilon_{n}(q)$ are determined by the orthogonality of the system $\Theta$ and are given by
$$\mu_{n}(q) = \frac{1}{\textrm{ce}_{n}(i\xi,q)\textrm{ce}_{n}(i\lambda,q)\textrm{ce}_{n}(\theta,q)}
\int_{-\pi}^{\pi}J_{0}(k\cdot\rho(\xi,\eta,\lambda,\theta))\textrm{ce}_{n}(\eta,q)d\eta, n\geq0,$$
$$\hskip-0.3cm\upsilon_{n}(q) = \frac{1}{\textrm{se}_{n}(i\xi,q)\textrm{se}_{n}(i\lambda,q)\textrm{se}_{n}(\theta,q)}
\int_{-\pi}^{\pi}J_{0}(k\cdot\rho(\xi,\eta,\lambda,\theta))\textrm{se}_{n}(\eta,q)d\eta, n\geq1$$
where we define $\upsilon_{0}(q) = 0$.

\section{The Reconstruction Problem and The Main Result}

Let $f$ be a continuous function, defined on $\Bbb R^{2}$, and assume that the SMT, $Rf$ of $f$, is given on the cylindrical surface $E\times\Bbb R^{+}$ where the ellipse $E$ is given by (1.1). Observe that by using rotation, translation and scaling we can all ways convert a general ellipse in the plane to be in this form. Hence, our data consists of all the integrals of $f$ on circles with centers on the ellipse $E$ and with arbitrary radii. Our main goal is to reconstruct $f$ from this data. Since the ellipse $E$ has the following parametrization
\vskip-0.2cm
$$\hskip-5cm E:(\cosh\xi_{0}\cos\theta, \sinh\xi_{0}\sin\theta), \theta\in[-\pi,\pi)$$
it follows that our data consists of the following set of integrals
$$\hskip-6.5cm (Rf)((\cosh\xi_{0}\cos\theta, \sinh\xi_{0}\sin\theta), r)$$
\begin{equation} \hskip-3cm = \int_{\Bbb S^{1}}f((\cosh\xi_{0}\cos\theta, \sinh\xi_{0}\sin\theta) + ry)rdS(y), \theta\in[-\pi,\pi), r\geq0\end{equation}
and our aim is to reconstruct $f$ from (2.1). For this, we have the following theorem.

\begin{thm}
Let $f$ be a continuous function with support inside the ellipse $E$. Let $x$ be an arbitrary point in $\Bbb R^{2}$ and let $x(\xi',\eta')$ be its presentation in elliptical coordinates, then
\vskip-0.2cm
$$\hskip-3cm f(x(\xi', \eta')) = \lim_{r'\rightarrow0}\frac{1}{2\pi}\int_{0}^{\infty}\Phi(\xi_{0}, \xi', \eta', q(k))J_{0}(r'k)kdk$$
where
\newpage
$$\hskip-10cm\Phi(\xi_{0}, \xi', \eta', q)$$
$$ = \sum_{n = 0}^{\infty}\left[\frac{\mathrm{ce}_{n}(i\xi',q)\mathrm{ce}_{n}(\eta',q)}{\mathrm{ce}_{n}(i\xi_{0}, q)}\int_{0}^{\infty}\int_{-\pi}^{\pi}(Rf)(x(\xi_{0},\eta), r)\mathrm{ce}_{n}(\eta,q)J_{0}(rk)d\eta dr \right.$$ $$\hskip1.5cm\left. + \frac{\mathrm{se}_{n}(i\xi',q)\mathrm{se}_{n}(\eta',q)}{\mathrm{se}_{n}(i\xi_{0}, q)}\int_{0}^{\infty}\int_{-\pi}^{\pi}(Rf)(x(\xi_{0},\eta), r)\mathrm{se}_{n}(\eta,q)J_{0}(rk)d\eta dr\right].$$
\end{thm}

For the proof of Theorem 2.1 we will follow the method introduced by Norton in \cite{22}. There, the author uses the following factorization
$$\hskip-1cm\delta(|re^{i\theta} - Re^{i\phi}| - \rho) = \rho\int_{0}^{\infty}J_{0}(k|re^{i\theta} - Re^{i\phi}|)J_{0}(k\rho)kdk$$
of the delta function and then the following expansion
\begin{equation}\hskip-2.3cm J_{0}(k|re^{i\theta} - Re^{i\phi}|) = \sum_{n = -\infty}^{\infty}J_{n}(kr)J_{n}(kR)e^{-in\theta}e^{in\phi}\end{equation}
of the Bessel function ($J_{n}$ denotes the Bessel function of order $n$) in order to find an explicit relation between each term in the spherical harmonic expansions of $f$ and $Rf$. The equation which defines the relation between each corresponding terms in the expansion can be inverted using the Hankel transform and thus each term in the spherical harmonic expansion of $f$ can be recovered. Thus, $f$ can also be recovered.

Here, we use the same procedure as described above where now, since we are working with elliptical coordinates, we will use the expansion (1.4) of the Bessel function $J_{0}$, which was proved by H.S. Cohl and H.Volkmer in \cite{8}, instead of the expansion (2.2) which was used by Norton in \cite{22}.

\section{Proof of the Main Result}

\textbf{Proof of Theorem 2.1}: From the definition of the spherical mean transform $Rf$ we have
\vskip-0.2cm
$$\hskip-6cm(Rf)((\cosh(\xi_{0})\cos(\eta), \sinh(\xi_{0})\sin(\eta)), r)$$
$$ \hskip-4cm = \int_{\Bbb S^{1}}f((\cosh(\xi_{0})\cos(\eta), \sinh(\xi_{0})\sin(\eta)) + ry)rdS(y)$$
$$ \hskip-3cm = \int_{\Bbb R^{2}}f(x)\delta\left(\left|(\cosh(\xi_{0})\cos(\eta), \sinh(\xi_{0})\sin(\eta)) - x\right| - r\right)dx.$$
In the last integral, making the following change of variables
$$\hskip-1.7cm x = (\cosh(\lambda)\cos(\theta), \sinh(\lambda)\sin(\theta)), dx = (\cosh^{2}\lambda - \cos^{2}\theta)d\theta d\lambda$$
to elliptical coordinates results in
\vskip-0.2cm
$$\hskip-6cm(Rf)((\cosh(\xi_{0})\cos(\eta), \sinh(\xi_{0})\sin(\eta)), r)$$
\begin{equation}
= \int_{0}^{\infty}\int_{-\pi}^{\pi}f(\cosh(\lambda)\cos(\theta),\sinh(\lambda)\sin(\theta))
\delta(\rho(\xi_{0},\eta, \lambda, \theta) - r)(\cosh^{2}\lambda - \cos^{2}\theta)d\theta d\lambda.
\end{equation}
where $\rho(\xi_{0}, \eta,\lambda, \theta) = |x(\xi_{0}, \eta) - x(\lambda, \theta)|$. Define the function $g$ by the relation
$$\hskip-2.5cm g(\lambda, \theta) = (\cosh^{2}\lambda - \cos^{2}\theta)f(\cosh(\lambda)\cos(\theta),\sinh(\lambda)\sin(\theta))$$
then using the following factorization
\vskip-0.2cm
$$\hskip-7cm \delta(\tau - r) = r\int_{0}^{\infty}J_{0}(\tau k)J_{0}(rk)kdk$$
of the delta function we have from equation (3.1) that
\vskip-0.2cm
$$\hskip-6cm(Rf)((\cosh(\xi_{0})\cos(\eta), \sinh(\xi_{0})\sin(\eta)), r)$$
\begin{equation}
 \hskip-3.2cm = r\int_{0}^{\infty}\int_{-\pi}^{\pi}\int_{0}^{\infty}g(\lambda, \theta)
J_{0}(k\cdot\rho(\xi_{0},\eta, \lambda, \theta))J_{0}(rk)kdk d\theta d\lambda.
\end{equation}
\begin{remark}
Observe that by our assumption on the support of $f$, the function $g$ is supported in the variable $\lambda$ and thus the last integral is well defined. For the rest of the proof we will use the fact that $g$ is supported with respect to the variable $\lambda$ and $Rf$ is supported with respect to the variable $r$ when integrating on infinite intervals and also when interchanging integration and summation. However, the reader should be aware that from now on we will not mention these facts each time that they are used.
\end{remark}
The right hand side of equation (3.2) is the Hankel transform of the function
$$\hskip-4cm\Lambda(k, \eta) = \int_{0}^{\infty}\int_{-\pi}^{\pi}g(\lambda, \theta)J_{0}(k\cdot\rho(\xi_{0},\eta, \lambda, \theta))d\theta d\lambda,$$
with respect to the variable $k$, evaluated at the point $r$. Thus, using the formula for the inverse of the Hankel transform we obtain that
$$\hskip-3.5cm\int_{0}^{\infty}\frac{1}{r}(Rf)((\cosh(\xi_{0})\cos(\eta), \sinh(\xi_{0})\sin(\eta)), r)J_{0}(kr)rdr$$
$$\hskip-6.5cm = \int_{0}^{\infty}\int_{-\pi}^{\pi}g(\lambda, \theta)J_{0}(k\cdot\rho(\xi_{0},\eta, \lambda, \theta))d\theta d\lambda$$
$$ \hskip-2.1cm = \frac{1}{\pi}\int_{0}^{\infty}\int_{-\pi}^{\pi}g(\lambda, \theta)
\sum_{n = 0}^{\infty}\left[\mu_{n}(q)\textrm{ce}_{n}(i\xi_{0}, q)\textrm{ce}_{n}(\eta,q)\textrm{ce}_{n}(i\lambda, q)\textrm{ce}_{n}(\theta,q)
\right.$$
$$\hskip-1.5cm\left. + \upsilon_{n}(q)\textrm{se}_{n}(i\xi_{0}, q)\textrm{se}_{n}(\eta,q)\textrm{se}_{n}(i\lambda, q)\textrm{se}_{n}(\theta,q)\right] d\theta d\lambda$$
$$ \hskip-1.1cm = \frac{1}{\pi}\sum_{n = 0}^{\infty}\left[\mu_{n}(q)\textrm{ce}_{n}(i\xi_{0}, q)\textrm{ce}_{n}(\eta,q)\int_{0}^{\infty}\int_{-\pi}^{\pi}g(\lambda, \theta)
\textrm{ce}_{n}(i\lambda, q)\textrm{ce}_{n}(\theta,q)d\theta d\lambda\right.$$
$$\left. + \upsilon_{n}(q)\textrm{se}_{n}(i\xi_{0}, q)\textrm{se}_{n}(\eta,q)\int_{0}^{\infty}\int_{-\pi}^{\pi}g(\lambda, \theta)
\textrm{se}_{n}(i\lambda, q)\textrm{se}_{n}(\theta,q)d\theta d\lambda\right]$$
where $k\in\Bbb R^{+}, \eta\in[-\pi,\pi)$ and where in the second passage we used the expansion (1.4) of the Bessel function $J_{0}$. Using the orthogonality of the system $\Theta$ and equation (1.3) yields for $n\geq0$:
\vskip-0.2cm
$$\hskip-7.7cm\int_{0}^{\infty}\int_{-\pi}^{\pi}g(\lambda, \theta)\textrm{ce}_{n}(i\lambda, q)\textrm{ce}_{n}(\theta, q)d\theta d\lambda$$
\begin{equation} = \frac{1}{\mu_{n}(q)\textrm{ce}_{n}(i\xi_{0}, q)}\int_{0}^{\infty}\int_{-\pi}^{\pi}(Rf)((\cosh(\xi_{0})\cos(\eta), \sinh(\xi_{0})\sin(\eta)), r)\textrm{ce}_{n}(\eta,q)J_{0}(rk)d\eta dr,\end{equation}
and for $n\geq1$
\vskip-0.2cm
$$\hskip-7.7cm\int_{0}^{\infty}\int_{-\pi}^{\pi}g(\lambda, \theta)\textrm{se}_{n}(i\lambda, q)\textrm{se}_{n}(\theta, q)d\theta d\lambda$$
\begin{equation} = \frac{1}{\upsilon_{n}(q)\textrm{se}_{n}(i\xi_{0}, q)}\int_{0}^{\infty}\int_{-\pi}^{\pi}(Rf)((\cosh(\xi_{0})\cos(\eta), \sinh(\xi_{0})\sin(\eta)), r)\textrm{se}_{n}(\eta,q)J_{0}(rk)d\eta dr.\end{equation}

Hence, for every $k\geq0$ (or equivalently $q\geq0$) and $n\geq0$ we acquired the following set of integrals
\vskip-0.2cm
$$\hskip-4.5cm\mathcal{K}_{n}(q) = \int_{0}^{\infty}\int_{-\pi}^{\pi}g(\lambda, \theta)\textrm{ce}_{n}(i\lambda, q)\textrm{ce}_{n}(\theta, q)d\theta d\lambda,$$
$$\hskip-4.5cm\mathcal{L}_{n}(q) = \int_{0}^{\infty}\int_{-\pi}^{\pi}g(\lambda, \theta)\textrm{se}_{n}(i\lambda, q)\textrm{se}_{n}(\theta, q)d\theta d\lambda.$$
Now, we take arbitrarily two series of functions $\{\rho_{n} = \rho_{n}(q)\}_{n = 0}^{\infty}$ and $\{\varrho_{n} = \varrho_{n}(q)\}_{n = 0}^{\infty}$
in the variable $q = \frac{1}{2}k^{2}$. Then, multiplying for each $n\geq0$ the integral $\mathcal{K}_{n}(q)$ with $\rho_{n}(q)$ and the integral $\mathcal{L}_{n}(q)$ with $\varrho_{n}(q)$ and then taking the sum we obtain the following expression
\vskip-0.2cm
$$\hskip-6.3cm\Psi(q) = \int_{0}^{\infty}\int_{-\pi}^{\pi}g(\lambda,\theta)K(\lambda,\theta,q)d\theta d\lambda$$
where
\begin{equation}K(\lambda,\theta,q) = \sum_{n = 0}^{\infty}\rho_{n}(q)\textrm{ce}_{n}(i\lambda, q)\textrm{ce}_{n}(\theta, q) + \varrho_{n}(q)\textrm{se}_{n}(i\lambda, q)\textrm{se}_{n}(\theta, q).\end{equation}

Now let $x$ be a fixed point in $\Bbb R^{2}$ and let $\xi'$ and $\eta'$ be the corresponding parameters in its parametrization in elliptical coordinates. Now, let us choose
$$\hskip-2cm \rho_{n}(q) = \mu_{n}(q)\textrm{ce}_{n}(i\xi',q)\textrm{ce}_{n}(\eta',q), \varrho_{n}(q) = \upsilon_{n}(q)\textrm{se}_{n}(i\xi',q)\textrm{se}_{n}(\eta',q)$$
and insert these expressions for $\rho_{n}(q)$ and $\varrho_{n}(q)$ in the expansion (3.5) of $K(\cdot,\cdot,q)$:

$$\hskip-4cm K(\lambda,\theta,q) = \sum_{n = 0}^{\infty}\left[\mu_{n}(q)\textrm{ce}_{n}(i\xi',q)\textrm{ce}_{n}(\eta',q)\textrm{ce}_{n}(i\lambda, q)\textrm{ce}_{n}(\theta, q)\right.$$ $$ \left.\hskip3cm + \upsilon_{n}(q)\textrm{se}_{n}(i\xi',q)\textrm{se}_{n}(\eta',q)\textrm{se}_{n}(i\lambda, q)\textrm{se}_{n}(\theta, q)\right] = J_{0}(k\cdot\rho(\xi',\eta', \lambda, \theta))$$
where in the last passage we used the expansion (1.4) for the Bessel function $J_{0}$. Hence, we obtained that
$$\hskip-1.5cm\int_{0}^{\infty}\int_{-\pi}^{\pi}g(\lambda,\theta)J_{0}(k\cdot\rho(\xi',\eta', \lambda, \theta))d\theta d\lambda = \int_{0}^{\infty}\int_{0}^{\pi}g(\lambda,\theta)K(\lambda,\theta,q)d\theta d\lambda$$
$$ \hskip-1.7cm = \sum_{n = 0}^{\infty}\left[\mu_{n}(q)\textrm{ce}_{n}(i\xi',q)\textrm{ce}_{n}(\eta',q)\int_{0}^{\infty}\int_{-\pi}^{\pi}g(\lambda,\theta)\textrm{ce}_{n}(i\lambda, q)\textrm{ce}_{n}(\theta, q)d\theta d\lambda \right.$$ $$\left.+ \upsilon_{n}(q)\textrm{se}_{n}(i\xi',q)\textrm{se}_{n}(\eta',q)\int_{0}^{\infty}\int_{-\pi}^{\pi}g(\lambda,\theta)\textrm{se}_{n}(i\lambda, q)\textrm{se}_{n}(\theta, q)d\theta d\lambda\right]$$
$$ \hskip-1.7cm = \sum_{n = 0}^{\infty}\left[\frac{\textrm{ce}_{n}(i\xi',q)\textrm{ce}_{n}(\eta',q)}{\textrm{ce}_{n}(i\xi_{0}, q)}\int_{0}^{\infty}\int_{-\pi}^{\pi}(Rf)(x(\xi_{0},\eta), r)\textrm{ce}_{n}(\eta,q)J_{0}(rk)d\eta dr \right.$$ \begin{equation}\hskip-0.5cm\left. + \frac{\textrm{se}_{n}(i\xi',q)\textrm{se}_{n}(\eta',q)}{\textrm{se}_{n}(i\xi_{0}, q)}\int_{0}^{\infty}\int_{-\pi}^{\pi}(Rf)(x(\xi_{0},\eta), r)\textrm{se}_{n}(\eta,q)J_{0}(rk)d\eta dr\right]\end{equation}
where in the last passage we used formulas (3.3) and (3.4). Let us denote by $\Phi = \Phi(\xi_{0}, \xi', \eta', q)$ the function defined by the infinite sum in the right hand side of equation (3.6). Multiplying equation (3.6) by $kr'J_{0}(r'k)$ (where $r' > 0$) and then integrating on $k\in[0,\infty)$ we have
$$ \int_{0}^{\infty}\Phi(\xi_{0}, \xi', \eta', q(k))J_{0}(r'k)r'kdk = \int_{0}^{\infty}\int_{0}^{\infty}\int_{-\pi}^{\pi}g(\lambda,\theta)J_{0}(k\cdot\rho(\xi',\eta', \lambda, \theta))J_{0}(kr')kr'd\theta d\lambda dk $$
$$ \hskip-4.7cm = \int_{0}^{\infty}\int_{-\pi}^{\pi}g(\lambda,\theta) \int_{0}^{\infty}J_{0}(k\cdot\rho(\xi',\eta', \lambda, \theta))J_{0}(kr')kr' dk d\theta d\lambda$$
$$ \hskip-7.2cm = \int_{0}^{\infty}\int_{-\pi}^{\pi}g(\lambda,\theta)\delta(\rho(\xi',\eta', \lambda, \theta) - r')d\theta d\lambda$$
$$\hskip-2.6cm = \left[y = (\cosh(\lambda)\cos(\theta), \sinh(\lambda)\sin(\theta)), dy = (\cosh^{2}\lambda - \cos^{2}\theta)d\theta d\lambda\right]$$
$$ \hskip-4.7cm = \int_{\Bbb R^{2}}f(y)\delta\left(\left|x(\xi', \eta') - y\right| - r'\right)dy = (Rf)\left(x(\xi', \eta'), r'\right).$$
Now, from the definition of the spherical mean transform it follows that
$$\hskip-7cm\lim_{r'\rightarrow0}\frac{1}{2\pi r'}(Rf)\left(x(\xi', \eta'), r'\right) = f(x(\xi',\eta')).$$
Since the point $x(\xi',\eta')$ is arbitrary we have proved Theorem 2.1. $\hskip3.2cm\square$

\end{document}